\documentclass[a4paper,11pt]{article}
\usepackage{amsfonts,amsmath,amssymb,amsthm,verbatim,placeins,xcolor,comment,placeins,caption,subcaption,enumitem,bm,graphicx,fancyvrb}
\usepackage{geometry}
\usepackage{float}
\usepackage[pdfauthor={DG},pdfstartview={FitH},colorlinks=true,citecolor=blue]{hyperref}
\usepackage[english]{babel}
\usepackage[numbers]{natbib}

\newtheorem{theorem}{Theorem}
\newtheorem{lemma}{Lemma}

\newtheorem{proposition}{Proposition}
\newtheorem{definition}{Definition}
\newtheorem{example}{Example}
\theoremstyle{remark}
\newtheorem{remark}{Remark}

\newcommand{\E}{\mathbb{E}}
\newcommand{\p}{\mathbb{P}}
\newcommand{\N}{\mathbb{N}}
\newcommand{\R}{\mathbb{R}}
\newcommand\1{\mathbf{1}}

\usepackage{booktabs}
\usepackage{bera}
\usepackage{xurl}
\usepackage{appendix}

\begin{document}
	
	\begin{center}
		\Large{\textbf{Integer-valued multifractal processes}}\footnote{
		This manuscript is a preprint of the paper: D.~Grahovac, Integer-valued multifractal processes, \textit{Chaos, Solitons \& Fractals} \textbf{201} (2025), no.~2, 117277, \url{https://doi.org/10.1016/j.chaos.2025.117277}\\
		© 2025. This manuscript version is made available under the CC-BY-NC-ND 4.0 license \url{https://creativecommons.org/licenses/by-nc-nd/4.0/}
		}\\
		\bigskip
		\bigskip
		Danijel Grahovac$^*$
	\end{center}
	
	\bigskip
	\begin{flushleft}
		\footnotesize{
			$^*$ School of Applied Mathematics and Informatics, J.~J.~Strossmayer University of Osijek, Trg Ljudevita Gaja 6, 31000 Osijek, Croatia; \texttt{dgrahova@mathos.hr}  [Corresponding author]
		}
		
	\end{flushleft}
	
	\bigskip
	
	\textbf{Abstract: } Multifractal scaling has been extensively studied for real-valued stochastic processes, but a systematic integer-valued analogue has remained largely unexplored. In this work, we introduce a multifractal framework for integer-valued processes using the thinning operation, which serves as a natural discrete counterpart to scalar multiplication. Within this framework, we construct integer-valued multifractal processes by time changing compound Poisson processes with nondecreasing multifractal clocks. We derive the scaling laws of their moments, provide explicit examples, and illustrate the results through numerical simulations. This construction integrates multifractal concepts into point process theory, enabling analysis of nonlinear discrete stochastic systems with nontrivial scaling properties.
	
	\medskip
	
	
	
	\medskip
	
	\noindent\textbf{Funding: } This work was supported by the Croatian Science Foundation under the project Scaling in Stochastic Models (HRZZ-IP-2022-10-8081).

	\bigskip

\section{Introduction} 

Multifractal processes have gained considerable attention over the past three decades since they have numerous applications in various fields \cite{mandelbrot1997mmar,abry2015irregularities,mandelbrot1997fractal,lovejoy2013weather,anh2009multifractal,jiang2019multifractal}. Multifractal scaling of moments can be used to analyze velocity or dissipation fields in the study of hydrodynamic turbulence \cite{frisch1985fully,robert2008hydrodynamic}. Multifractal models have successfully explained many phenomena in climatology \cite{lovejoy2013weather,kalamaras2017multifractal,laib2018multifractal,BARANOWSKI2019318}. In medicine, multifractal properties are observed in many dynamic biological processes \cite{pavlov2016multifractality,joseph2021multifractal}. Multifractal methods proved very useful in image analysis \cite{abry2015irregularities,martsepp2022dependence,lopes2009fractal}. Financial data typically exhibits multifractal features, hence multifractal processes play a prominent role in finance \cite{mandelbrot1997mmar,bacry2008continuous,duchon2012forecasting,ZHOU2012147}. In engineering, multifractal analysis has proven effective across diverse applied contexts, like e.g.~communication systems \cite{gadre2003multifractal,wu2015low} or structural engineering \cite{ebrahimkhanlou2016multifractal} (see also Ref.~\cite{vehel2012fractals} for a survey). In physics, multifractality remains a powerful tool for characterizing complex, nonlinear dynamics \cite{apolinario2022dynamical}. Nonlinearity is crucial because it allows systems to exhibit behaviors such as chaos, self-organization, and emergent patterns (see Refs.~\cite{strogatz2024nonlinear,gao2025inhomogeneity,feng2025bilinear,liu2025n,WANG2026109720} for such examples).

There are many different ways on how to define a multifractal process \cite{riedi2003multifractal,muzy2002multifractal,grahovac2018bounds}. Regardless of the approach, multifractality is viewed as a scale invariance property that extends the concept of self-similarity, whether it be in terms of scaling of finite-dimensional distributions, or moments of the process (as seen, e.g.~in Refs.~\cite{mandelbrot1997mmar}, \cite{riedi2003multifractal}).

In terms of finite-dimensional distributions, a stochastic process $X=\{X(t), \, t \geq 0\}$ is said to be \textbf{multifractal} if for every $\lambda \in (0,1)$ there exists a positive random variable $M(\lambda)$, independent of $X$ such that
\begin{equation}\label{eq:mfdef}
	\{ X(\lambda t) \}_{t \in (0,T]} \overset{d}{=} \{ M(\lambda) X(t) \}_{t \in (0,T]}.
\end{equation}
Here, and in what follows, $\{\cdot\} =^d \{\cdot\}$ stands for the equality of finite-dimensional distributions of two stochastic processes. The parameter $T$ limits the time set over which the scaling property holds. The definition can be made even more general by taking an arbitrary set of scales $\Lambda$ instead of $(0,1)$ and an arbitrary set $\mathcal{S}$ instead of $(0,T]$ on which the scaling property holds. However, as was shown in Ref.~\cite{grahovac2020}, under very mild assumptions the general definition reduces to the typical case defined above. If $M(\lambda)$ is deterministic, then one necessarily has $M(\lambda)=\lambda^H$ for some $H\geq 0$ and the scaling property \eqref{eq:mfdef} reduces to \textbf{self-similarity}.

An important example of a truly multifractal process satisfying scale invariance \eqref{eq:mfdef} is the multiplicative cascade \cite{barral2002multifractal,bacry2003log,muzy2002multifractal}. Mandelbrot introduced cascades using a discrete grid-based construction \cite{mandelbrot1972}. Various similar constructions have been suggested to achieve continuous scaling properties, starting with Ref.~\cite{barral2002multifractal}, followed by Refs.~\cite{bacry2003log}, \cite{muzy2002multifractal}, and more recently, Ref.~\cite{barral2014exact}. 

Time series related to counts occur across a variety of applications, ranging from medical science, epidemiology, meteorology, network modeling to actuarial science, econometrics, and finance (see, e.g.~Refs.~\cite{davis1998,cameron_trivedi_1998,weiss2008,mckenzie2003,leonenko2007,Hilbe_2014} and the references therein). The widespread presence of multifractal phenomena indicates that multifractal models could be quite valuable for count data. 

The goal of this paper is to define multifractality for integer-valued processes. To our knowledge, no systematic study of integer-valued multifractal models has ever been conducted. Some models in this context were proposed in unpublished manuscripts \cite{boucher2013time,boucher2015time}, where a Poisson process model is considered with random intensity involving multifractal process. A similar approach has been considered in spatial setting in Ref.~\cite{baile2021}. However, a general theoretical framework for multifractality in integer-valued processes is still missing, and research in this direction remains largely unexplored. This gap motivates the present work, where we aim to formalize the concept and provide rigorous foundations for future applications. 

Clearly, \eqref{eq:mfdef} cannot hold for integer-valued processes unless $M(\lambda)$ is integer-valued for all $\lambda\in (0,1)$. Moreover, cascade processes, being the main example of multifractal processes, do not have atoms (see Lemma 1.4 in Ref.~\cite{barral2014exact}). To define multifractality for integer-valued processes, we need to use the multiplication operation preserving integer values instead of the classical multiplication. Our approach is based on the thinning approach introduced in Ref.~\cite{steutel1979} and generalized thinning as defined in Ref.~\cite{vaHarnVervaat1982}, which we further extend to incorporate random multiplicative factors. Instead of the distributional approach, we use the connection with continuous-time branching process to extend the definition to stochastic processes.

We next construct a broad class of integer-valued multifractal processes by applying a random time-change to a compound Poisson process, where the new random clock is a process that is multifractal in the classical sense. Within this framework, we derive scaling relations for the moments of the resulting processes and illustrate the results through concrete examples and numerical simulations.

The models developed in this paper open a possibility for modeling count data, which arises in diverse contexts such as epidemiology, finance, biology, and physics \cite{davis1998,cameron_trivedi_1998,Hilbe_2014}. By introducing multifractality into integer-valued processes, our work addresses a gap in theory while providing tools for applications where discrete data with complex variability are central.

\section{Thinning operators}
Before discussing different approaches to the multiplication operation that preserves integer values, we introduce some notation. In what follows, $\N_0$ stands for the set of nonnegative integers and $\R_+=[0,\infty)$ for nonnegative reals. For an $\N_0$-valued random variable $X$, $P_X$ denotes its probability generating function (\textbf{pgf})
\begin{equation*}
	P_X(z)=\sum_{k=0}^\infty z^k \p (X=k), \quad |z|\leq 1.
\end{equation*}

\subsection{Steutel \& van Harn operators}
Since the multiplication of an integer by a real number can result in a noninteger value, in Ref.~\cite{steutel1979} a multiplication operation is introduced to replace scalar multiplication. For $\alpha \in (0,1)$ and an $\N_0$-valued random variable $X$, the \textbf{Steutel \& van Harn operation} $\odot$ is defined in distribution by
\begin{equation}\label{eq:bin-thin}
	P_{\alpha \odot X}(z) = P_X(1 - \alpha + \alpha z).
\end{equation}
It follows that 
\begin{equation*}
	\alpha \odot X \overset{d}{=} Z_1 + \cdots + Z_X,
\end{equation*}
where $Z_1,Z_2,\dots$ are independent random variables, independent of X, all having Bernoulli distribution with parameter $\alpha$, i.e.~$\p(Z_i=1)=1-\p(Z_i=0)=\alpha$. Conditionally on $X$, $\alpha \odot X$ has a binomial distribution with parameters $X$ and $\alpha$. Therefore, $\odot$ is often called a \textbf{binomial thinning operation}.

The multiplication operation \eqref{eq:bin-thin} has been generalized in Ref.~\cite{vaHarnVervaat1982}. To this end, let $F=(F_s)_{s\geq0}$ be a continuous composition semigroup of pgfs, i.e.~a collection of non-constant pgfs such that
\begin{equation*}
	F_s \circ F_t = F_{s+t}, \quad s,t \geq 0,
\end{equation*}
and $\lim_{s\to 0} F_s = F_0=I$, where $I(z)=z$, $|z|\leq 1$. We also assume that 
\begin{equation}\label{e:Fassum}
	F_t\neq I, t>0 \text{ and } F_1'(1)=e^{-1},
\end{equation}
which implies that $F_s'(1)=e^{-s}$, $s > 0$. For $\alpha \in (0,1)$ and $\N_0$-valued random variable $X$, the \textbf{generalized Steutel \& van Harn operator} induced by the semigroup $F$, denoted by $\odot_F$, is defined in distribution by
\begin{equation}\label{eq:gen-thin}
	P_{\alpha \odot_F X} (z) = P_X (F_{s}(z)), \quad s=-\log \alpha.
\end{equation}
For the semigroup $F_s(z)=1-e^{-s}+e^{-s}z$, we get the standard multiplication operation \eqref{eq:bin-thin}.

The probabilistic interpretation relies on the connection of the semigroups with branching processes (see e.g.~Refs.~\cite{harris1963,AthreyaNey1972}). Every semigroup $F$ determines a continuous-time Markov branching process $Z=\{Z(t),\, t \geq 0\}$ such that $Z(0)=1$ a.s.~and
\begin{equation*}
	P_{Z(s)} (z) = F_s(z), \quad s \geq 0.
\end{equation*}
We call $Z$ an \textbf{$F$-branching process}. The assumption $F_1'(1)=e^{-1}$ implies that $\E Z(1)<1$, so the branching process is subcritical and $\lim_{s\to \infty} F_s(z) = 1$, $z\geq 0$. We can now interpret \eqref{eq:gen-thin} as
\begin{equation*}
	e^{-s} \odot_F X \overset{d}{=} Z_1(s) + Z_2(s) + \cdots + Z_X(s) =: T_X(s),
\end{equation*}
where $Z_1, Z_2,\dots$ are independent $F$-branching processes, independent of $X$. The process $T_X$ is a branching process such that $T_X(0)=X$.

\medskip

The previous definition can be extended to include a random multiplier $A$ with values in $(0,1)$, independent of $X$, replacing the deterministic $\alpha$ (see e.g.~Ref.~\cite[Appendix A]{steutelvanharn2004book}). In this case, the distribution of the multiplication operator is simply a mixture of the form
\begin{equation}\label{eq:gen-thin-rand}
	P_{A \odot_F X} (z) = \int_{(0,1)} P_X (F_{-\log \alpha} (z)) \p_A(d \alpha),
\end{equation}
where $\p_A$ is the distribution of $A$.

\subsection{Thinning for stochastic processes}
To adapt the definition of multifractality for integer-valued processes, we need to specify how the thinning operators act on the higher-dimensional distributions of the process. For standard multiplication, it is clear that
\begin{equation*}
	\phi_{\alpha X(t_1), \alpha X(t_2)}(z_1,z_2) = \phi_{X(t_1),X(t_2)}(\alpha z_1, \alpha z_2),
\end{equation*}
where $\phi_{Y_1,Y_2}$ denotes the characteristic function of the random vector $(Y_1,Y_2)$. Taking the same approach would mean putting
\begin{equation}\label{e:wrongapproach}
	P_{\alpha \odot_F X(t_1), \alpha \odot_F X(t_2)}(z_1,z_2) = P_{X(t_1),X(t_2)}(F_s (z_1), F_s (z_2)), \quad s=-\log \alpha,
\end{equation}
where $P_{Y_1,Y_2}$ denotes bivariate pgf. It turns out that this approach leads to a somewhat pathological class of processes. Namely, as elaborated in Ref.~\cite{vanharn1985selfsimilar}, even for integer-valued self-similar processes this could lead to a non-measurable process. To avoid this issue, the solution is to define the process obtained after multiplication in a pathwise way and not only distributionally, as it was done in \eqref{eq:bin-thin} and \eqref{eq:gen-thin} for random variables.

\begin{definition}\label{def:processthin1}
	Let $X=\{X(t), \, t \geq 0\}$ be an $\N_0$-valued stochastic process and for some continuous semigroup $F=(F_s)_{s\geq 0}$ satisfying \eqref{e:Fassum} suppose that $Z_1, Z_2,\dots$ are independent $F$-branching processes, independent of $X$. Then for $\alpha \in (0,1)$ we define the process $\{\alpha \odot_F X(t), \, t \geq 0\}$ by
	\begin{equation*}
		\alpha \odot_F X(t) = T_{X(t)} (- \log \alpha), \quad t \geq 0,
	\end{equation*}
	where
	\begin{equation*}
		T_{X(t)}(s) = Z_1(s) + \cdots + Z_{X(t)}(s).
	\end{equation*}
\end{definition}

This definition is clearly consistent with \eqref{eq:gen-thin}. However, the bivariate distributions 
\begin{equation*}
	(\alpha \odot_F X(t_1), \alpha \odot_F X(t_2))=(T_{X(t_1)}(-\log \alpha), T_{X(t_2)}(-\log \alpha) )
\end{equation*}
do not have simple expression as in \eqref{e:wrongapproach} since we use the same branching process for both $X(t_1)$ and $X(t_2)$.

The previous approach was used in Ref.~\cite{vanharn1985selfsimilar} to define integer-valued self-similar processes. We generalize Definition \ref{def:processthin1} to allow for a random multiplier.

\begin{definition}\label{def:processthin-rand}
	Let $X=\{X(t), \, t \geq 0\}$ be an $\N_0$-valued stochastic process and for some continuous semigroup $F=(F_s)_{s\geq 0}$ satisfying \eqref{e:Fassum} suppose that $Z_1, Z_2,\dots$ are independent $F$-branching processes, independent of $X$. Then for $(0,1)$-valued random variable $A$ independent of $X$ and $Z_1, Z_2,\dots$, we define the process $\{A \odot_F X(t), \, t \geq 0\}$ by
	\begin{equation*}
		A \odot_F X(t) = T_{X(t)} (- \log A), \quad t \geq 0,
	\end{equation*}
	where
	\begin{equation}\label{eq:Tdef}
		T_{X(t)}(s) = Z_1(s) + \cdots + Z_{X(t)}(s).
	\end{equation}
\end{definition}

The previous definition is in line with \eqref{eq:gen-thin-rand} since by conditioning on $A$ we get
\begin{equation*}
	P_{A \odot_F X(t)}(z) = \int_{(0,1)} P_{X(t)} (F_{-\log \alpha} (z)) \p_A(d \alpha),
\end{equation*}
where $\p_A$ denotes the distribution of $A$.

Note that for each fixed process $X$ we take one sequence of branching processes in the multiplication operation. In repeated application of this operation, we then consider independent sequences for each operation. With this interpretation, we have the following. 

\begin{proposition}\label{prop:ABrepeated}
	Let $A$ and $B$ be independent $(0,1)$-valued random variables, independent of $\N_0$-valued stochastic process $X=\{X(t), \, t \geq 0\}$. Then
	\begin{equation*}
		\{ A \odot_F (B \odot_F X(t))\} \overset{d}{=} \{(A B) \odot_F X(t)\}.
	\end{equation*}
\end{proposition}

\begin{proof}
	Let $Z'_1, Z'_2, \dots$ and $Z''_1, Z''_2, \dots$ denote the $F$-branching processes used in the definition of $B \odot_F X(t)$ and $A \odot_F (B \odot_F X(t))$, respectively. Furthermore, let $T'_n(s)=Z'_1(s)+\cdots+Z'_n(s)$ and $T''_n(s)=Z''_1(s)+\cdots+Z''_n(s)$. Then we have that
	\begin{equation*}
		A \odot_F (B \odot_F X(t)) = Z''_1(-\log A)+\cdots+Z''_{T'_{X(t)}(-\log B)}(-\log A),
	\end{equation*}
	which can be written as
	\begin{align*}
		A \odot_F & (B \odot_F X(t)) = Z''_1(-\log A)+\cdots+Z''_{T'_{1}(-\log B)}(-\log A)\\
		&+ Z''_{T'_{2}(-\log B)+1}(-\log A)+\cdots+Z''_{T'_{2}(-\log B)}(-\log A) + \cdots\\
		&+ Z''_{T'_{X(t)-1}(-\log B)+1}(-\log A)+\cdots+Z''_{T'_{X(t)}(-\log B)}(-\log A).
	\end{align*}
	If we denote
	\begin{align*}
		&Z_j(s) =\\
		&\ \begin{cases}
			Z_j'(s), & \text{ if } s \leq - \log B,\\
			Z''_{T'_{j-1}(-\log B)+1}(s + \log B) + \dots + Z''_{T'_j(-\log B)}(s + \log B),& \text{ if } s > - \log B,\\
		\end{cases}
	\end{align*}
	then
	\begin{equation}\label{eq:proof:mult}
		\begin{aligned}
			A \odot_F (B \odot_F X(t)) =& Z_1(-\log A-\log B) + Z_2(-\log A - \log B)\\
			& + \cdots + Z_{X(t)}(-\log A - \log B).
		\end{aligned}
	\end{equation}
	Due to the branching property, the processes $Z_1, Z_2,\dots$ are independent $F$-branching processes, independent of $X$. On the other hand
	\begin{equation*}
		(AB) \odot_F X(t) =  Z'_1(-\log A - \log B) + \cdots + Z'_{X(t)}(-\log A - \log B),
	\end{equation*}
	and this process has the same finite-dimensional distributions as \eqref{eq:proof:mult}.
\end{proof}

\begin{remark}
	The equality of one-dimensional distributions in Proposition \ref{prop:ABrepeated} follows easily by the semigroup property since
	\begin{align*}
		P_{A \odot_F (B \odot_F X(t))} (z) &= \int_{(0,1)} P_{B \odot_F X(t)} (F_{-\log \alpha} (z)) \p_A(d \alpha)\\
		&=\int_{(0,1)} \int_{(0,1)} P_{X(t)} (F_{-\log \beta} (F_{-\log \alpha} (z))) \p_A(d \alpha) \p_B(d \beta)\\        
		&=\int_{(0,1)} \int_{(0,1)} P_{X(t)} (F_{-\log (\alpha \beta)} (z)) \p_A(d \alpha) \p_B(d \beta)\\
		&=P_{(A B) \odot_F X(t)}(z).
	\end{align*}
\end{remark}

\section{Definition of integer-valued multifractal processes}

For the definition of integer-valued multifractal process we use multiplication operation from Definition \ref{def:processthin-rand} as follows.

\begin{definition}
	Let $F=(F_s)_{s\geq 0}$ be a continuous semigroup satisfying \eqref{e:Fassum}. An $\N_0$-valued stochastic process $X=\{X(t), \, t \geq 0\}$ is \textbf{$F$-multifractal} if for some $T>0$ and every $\lambda\in (0,1)$ there exists a $(0,1)$-valued random variable $M(\lambda)$, independent of $X$, such that
	\begin{equation*}
		\{ X(\lambda t) \}_{t \in (0,T]} \overset{d}{=} \{ M(\lambda) \odot_F X(t) \}_{t \in (0,T]}.
	\end{equation*}
\end{definition}

If $M(\lambda)=m(\lambda)$ is deterministic for every $\lambda\in (0,1)$, then the previous definition reduces to the definition of integer-valued \textbf{$F$-self-similar} process given in Ref.~\cite{vanharn1985selfsimilar}. In this case, the scaling property holds over the time set $(0,\infty)$ and if $X$ is additionally continuous in distribution, then $m(\lambda)=\lambda^H$ for some $H\geq 0$.

The definition of classical multifractal process \eqref{eq:mfdef} can be made more general by allowing the scaling factors to additionally depend on $t$, as was done in Ref.~\cite{grahovac2020}. Namely, one can postulate that for every $\lambda \in (0,1]$ there exists a family of identically distributed positive random variables $\{M(\lambda,t), \, t \in (0,T] \}$, independent of $X$ such that
\begin{equation}\label{eq:mfdef-gen}
	\{ X(\lambda t) \}_{t \in (0,T]} \overset{d}{=} \{ M(\lambda, t) X(t) \}_{t \in (0,T]}.
\end{equation}
Adapting this definition would require defining integer-valued multiplication for two processes. For a $(0,1)$-valued process $A=\{A(t), \, t\geq 0\}$, independent of $X$, we could do this by putting, similarly as in Definition \ref{def:processthin-rand},
\begin{equation*}
	A(t) \odot_F X(t) = T_{X(t)} (- \log A(t)), \quad t \geq 0.
\end{equation*}
We do not pursue this approach here, since for all the main examples of multifractals the family of scaling factors is actually independent of $t$. In particular, for cascade processes we have $M(\lambda, t)=M(\lambda)$ for any $t\in (0,T]$.

\begin{remark}\label{rem:directcosntr}
	We note that it is possible to construct an integer-valued process that satisfies \eqref{eq:mfdef-gen} by using the approach described in Ref.~\cite{grahovac2020}. In fact, let $L=\{L(t), \, t \geq 0\}$ be a compound Poisson process with jump distribution having values in the set $K=\{\log k \colon k=1,2,\dots \}$. In particular, $L$ is a $K$-valued L\'evy process. By Theorem 1 in Ref.~\cite{grahovac2020}, the process $X^{(a)}(t)=\exp \{L(a+\log t) - L(a)\}$, $t \in (0, e^a]$, can be extended to a process $X=\{X(t), \, t\geq 0\}$ satisfying \eqref{eq:mfdef-gen} such that $\{X(t), \, (0, e^a]\} \overset{d}{=} \{X^{(a)}(t), \, (0, e^a]\}$, for any $a\geq 0$. Since $L$ is $K$-valued, the obtained process $X$ will be $\N_0$-valued.
\end{remark}

Multifractality is frequently defined as the scaling property of the moments of the process. More precisely, the multifractal scaling of moments holds if for some range of values of $q$ and $t\in (0,T]$ it holds that
\begin{equation}\label{eq:momscal}
	\E|X(t)|^q = c(q) t^{\tau (q)},
\end{equation}
where the so-called scaling function $q \mapsto \tau(q)$ is nonlinear, as opposed to self-similar processes for which the scaling function is always linear $\tau(q)=Hq$. Multifractal processes satisfying \eqref{eq:mfdef} always have an exact scaling of moments as in \eqref{eq:momscal} and in all known examples \eqref{eq:momscal} is a consequence of \eqref{eq:mfdef}. We study scaling of moments in the next section.

\section{Construction}

Through this section we will fix a semigroup $F=(F_s)_{s\geq 0}$ satisfying \eqref{e:Fassum}, related to a branching process $Z=\{Z(t), \, t \geq 0\}$. Since the branching process is subcritical, we have by Theorem IV.2.4 in Ref.~\cite{asmussenhering1983} that the limits
\begin{equation}\label{eq:Gdef}
	G(z) = \lim_{t\to \infty} \frac{F_t(z) - F_t(0)}{1-F_t(0)}, \quad 0\leq z \leq 1,
\end{equation}
exist and $G$ is a pgf that is a unique solution of 
\begin{equation}\label{eq:Gproperty}
	1-G(F_t(z)) = e^{-t} (1-G(z)), \quad G(0)=0,
\end{equation}
for any $t>0$.

For the construction of integer-valued multifractal process, we take a specific compound Poisson process and make a time change with some process that is multifractal in the classical sense. To this end, let $N=\{N(t), \, t \geq 0\}$ denote the compound Poisson process with rate parameter $1$ and jump distribution determined by the pgf $G$ given in \eqref{eq:Gdef}. This means that $N$ has the form
\begin{equation}\label{eq:defPP}
	N(t) = W_1 + \cdots + W_{N_1(t)}, \quad t \geq 0,
\end{equation}
where $N_1$ is a unit rate Poisson process and $W_1,W_2,\dots$ are independent random variables with pgf $G$. In particular, we have
\begin{equation*}
	P_{N(t)}(z) = \exp \{-t (1-G(z))\}.
\end{equation*}
By Lemma 4.1 in Ref.~\cite{vanharn1985selfsimilar}, the finite-dimensional distributions of $N$ are given by
\begin{equation}\label{eq:Nfidispgf}
	\begin{aligned}
		&P_{N(t_1),\dots, N(t_n)} (z_1, \dots, z_n) = \E \left( z_1^{N(t_1)} \cdots z_n^{N(t_n)}\right)\\
		&\ = \exp \big\{ -t_1 (1-G(z_1\cdots z_n)) - (t_2-t_1) (1-G(z_2\cdots z_n))\\
		&\quad \quad \quad \quad - \cdots - (t_n-t_{n-1}) (1-G(z_n))\big\},
	\end{aligned}
\end{equation}
for $n\in \N$ and $0\leq t_1 \leq \cdots \leq t_n$. Moreover, $N$ is $F$-self-similar with $H=1$, i.e.
\begin{equation}\label{eq:Nss-prop}
	\left\{ N(\lambda t) \right\} \overset{d}{=} \left\{ \lambda  \odot_F N(t) \right\}.
\end{equation}
We now show that this property holds for random multipliers, too. First, we need a technical lemma.

\begin{lemma}\label{lemma:conditioning}
	Let $\bm{U}=\{\bm{U}(t), \, t \geq 0\}$, $\bm{U}(t)=(U_{1}(t),\dots,U_{n}(t))$, be an $\R^n$-valued process with right-continuous sample paths and $V_1,\dots,V_n$ random variables, independent of $\bm{U}$. Then for any bounded continuous function $f\colon \R^n \to \R$ it holds
	\begin{equation*}
		\E f(U_{1}(V_1),\dots,U_{n}(V_n)) = \E g(V_1,\dots,V_n),
	\end{equation*}
	where $g(v_1,\dots,v_n)=\E f(U_{1}(v_1),\dots,U_{n}(v_n))$.
\end{lemma}

\begin{proof}
	First note that, since $\bm{U}(t,\omega)$ is measurable in $\omega$ and right-continuous in $t$, it is jointly measurable, i.e.~$(t, \omega) \mapsto \bm{U}(t,\omega)$ is $\mathcal{B}(\R_+) \times \mathcal{F}$-measurable, where $\mathcal{F}$ is the $\sigma$-algebra of the underlying probability space. Hence, $(U_{1}(V_1),\allowbreak\dots,U_{n}(V_n))$ is a well-defined random vector. 
	
	Let $h_k(s)=\sum_{j=1}^{k^2} (j/k) \1_{[(j-1)/k, j/k)}(s)$. Then $f(U_{1}(h_k(V_1)),\dots,\allowbreak U_{n}(h_k(V_n)))$ is a measurable function of $U_{i}(1/k), U_{i}(2/k), \dots, U_{i}(k)$, $i=1,\dots,n$ and $V_1,\dots,V_n$. We can apply Proposition 1.16 from Ref.~\cite{sato1999levy} to conclude that
	\begin{equation*}
		\E f(U_{1}(h_k(V_1)),\dots,U_{n}(h_k(V_n))) = \E g_k(V_1,\dots,V_n),
	\end{equation*}
	where $g_k(v_1,\dots,v_n)= \E f(U_{1}(h_k(v_1)),\dots,U_{n}(h_k(v_n)))$. Since $h_k(V_i) \downarrow V_i$, $i=1,\dots,n$, as $k\to \infty$, by the right continuity of paths of $\bm{U}$ we have
	\begin{equation}\label{eq:lemma:cond:limit}
		\begin{aligned}
			\E f(U_{1}(V_1),\dots,U_{n}(V_n)) &= \lim_{k\to \infty} \E f(U_{1}(h_k(V_1)),\dots,U_{n}(h_k(V_n)))\\
			&= \lim_{k\to \infty} \E g_k(V_1,\dots,V_n) = \E g(V_1,\dots,V_n).
		\end{aligned}
	\end{equation}
\end{proof}

\begin{lemma}\label{lemma:Nscaling}
	The compound Poisson process $N=\{N(t), \, t \geq 0\}$ given by \eqref{eq:defPP} satisfies
	\begin{equation*}
		\left\{ N(A t) \right\} \overset{d}{=} \left\{ A  \odot_F N(t) \right\},
	\end{equation*}
	for any $(0,1)$-valued random variable $A$ independent of $N$.
\end{lemma}

\begin{proof}
	Let $Z_1, Z_2,\dots$ denote independent $F$-branching processes, independent of $X$ and $A$, defining the multiplication $\odot_F$. Then
	\begin{equation*}
		A \odot_F N(t) = Z_1(-\log A) + \cdots + Z_{N(t)}(-\log A) = T_{N(t)}(-\log A),
	\end{equation*}
	for $T$ given by \eqref{eq:Tdef}. We first prove \eqref{eq:Nss-prop} for completeness. For $0\leq t_1 \leq \cdots \leq t_n$
	\begin{align*}
		&P_{a \odot_F N(t_1),\dots, a \odot_F N(t_n)} (z_1, \dots, z_n) = \E \left( z_1^{a \odot_F N(t_1)} \cdots z_n^{a \odot_F N(t_n)}\right)\\
		&\ = \E \left( z_1^{T_{N(t_1)}(-\log a)} \cdots z_n^{T_{N(t_n)}(-\log a)}\right)\\
		&\ = \E \Big( (z_1\cdots z_n)^{Z_1(-\log a) + \cdots + Z_{N(t_1)}(-\log a)} (z_2\cdots z_n)^{Z_{N(t_1)+1}(-\log a) + \cdots + Z_{N(t_2)}(-\log a)}\\
		&\hspace{2cm}  \cdots z_n^{Z_{N(t_{n-1})+1}(-\log a) + \cdots + Z_{N(t_n)}(-\log a)} \Big)\\
		&\ = P_{N(t_1)}(F_{-\log a}(z_1 \cdots z_n)) P_{N(t_2)-N(t_1)}(F_{-\log a}(z_2 \cdots z_n))\\
		&\hspace{2cm} \cdots P_{N(t_n)-N(t_{n-1})}(F_{-\log a}(z_n))\\
		&\ = \exp \Big\{ -t_1 (1-G(F_{-\log a}(z_1 \cdots z_n))) - (t_2-t_1) (1-G(F_{-\log a}(z_2 \cdots z_n)))\\
		&\hspace{2cm}  - \cdots - (t_n-t_{n-1}) (1-G(F_{-\log a}(z_n)))\Big\}\\
		&\ = \exp \big\{ -a t_1 (1-G(z_1\cdots z_n)) - a(t_2-t_1) (1-G(z_2\cdots z_n))\\
		&\hspace{2cm} - \cdots - a(t_n-t_{n-1}) (1-G(z_n))\big\}\\
		&\ = P_{N(at_1),\dots, N(at_n)} (z_1, \dots, z_n)
	\end{align*}
	where the last two lines follow by \eqref{eq:Gproperty} and \eqref{eq:Nfidispgf}. Now  
	\begin{align*}
		P_{A \odot_F N(t_1),\dots, A \odot_F N(t_n)} (z_1, \dots, z_n) &= \E \left( z_1^{T_{N(t_1)}(-\log A)} \cdots z_n^{T_{N(t_n)}(-\log A)}\right)\\
		&=\E f(\bm{U}(A)),
	\end{align*}
	where $f(u_1,\dots,u_n)=z_1^{u_1}\cdots z_n^{u_n}$ and $\bm{U}(v)=(T_{N(t_1)}(-\log v),\allowbreak \dots, \allowbreak T_{N(t_n)}(-\log v))$. By Lemma \ref{lemma:conditioning} applied on $\bm{U}$ and $V_1=\cdots=V_n=A$, we get
	\begin{equation*}
		P_{A \odot_F N(t_1),\dots, A \odot_F N(t_n)} (z_1, \dots, z_n) = \E g(A),
	\end{equation*}
	for 
	\begin{equation*}
		g(a)=P_{a \odot_F N(t_1),\dots, a \odot_F N(t_n)} (z_1, \dots, z_n) = P_{N(at_1),\dots, N(at_n)} (z_1, \dots, z_n).
	\end{equation*}
	On the other hand, for $\bm{U}(u)=(N(u t_1),\dots, N(u t_n))$ and $V_1=\cdots=V_n=A$, Lemma \ref{lemma:conditioning} gives
	\begin{equation*}
		P_{N(At_1),\dots, N(At_n)} (z_1, \dots, z_n) = \E f(\bm{U}(A)) = \E g(A),
	\end{equation*}
	which proves the statement.
\end{proof}


To construct integer-valued multifractals we make a composition of the compound Poisson process $N$ and some classically multifractal nondecreasing process $Y$. Note that since $Y$ is nondecreasing, $M(\lambda)$ in \eqref{eq:mfdef} is a.s.~$(0,1)$-valued and the integer-valued multiplication with such factors will be well-defined.

\begin{theorem}\label{thm:constr}
	Let $Y=\{Y(t), \, t \geq0\}$ be an $\R_+$-valued nondecreasing multifractal process satisfying \eqref{eq:mfdef} and $N=\{N(t), \, t \geq 0\}$ a compound Poisson process with rate $1$ and jump distribution $G$, independent of $Y$. Then an $\N_0$-valued  process $X=\{X(t), \, t \geq 0\}$ defined by
	\begin{equation*}
		X(t) = N(Y(t)), \quad t \geq 0,
	\end{equation*}
	is $F$-multifractal with $F$ and $G$ related by \eqref{eq:Gproperty}.
\end{theorem}

\begin{proof}
	First, note that we can always take the right-continuous modification of $N$. Then $N$ is jointly measurable and therefore $X$ is well-defined. By the scaling property of $Y$ we have
	\begin{equation*}
		\{ X(\lambda t) \} = \{ N(Y(\lambda t)) \} \overset{d}{=} \{ N(M(\lambda) Y(t)) \},
	\end{equation*}
	where $N$, $M(\lambda)$ and $Y$ are independent. Now, for arbitrary $0\leq t_1 \leq \cdots \leq t_n$ and $z_1,\dots,z_n$,
	\begin{equation*}
		P_{N(M(\lambda) Y(t_1)),\dots, N(M(\lambda) Y(t_n))} (z_1, \dots, z_n) = \E f(U_1(Y(t_1)), \dots, U_n(Y(t_n))),
	\end{equation*}
	for $f(u_1,\dots,u_n)=z_1^{u_1}\cdots z_n^{u_n}$ and $U_i(u)=N(M(\lambda) u)$, $i=1,\dots,n$. By Lemma \ref{lemma:conditioning} for $V_i=Y(t_i)$, $i=1,\dots,n$, we get
	\begin{equation*}
		P_{N(M(\lambda) Y(t_1)),\dots, N(M(\lambda) Y(t_n))} (z_1, \dots, z_n) = \E g(Y(t_1), \dots, Y(t_n)),
	\end{equation*}
	where, by Lemma \ref{lemma:Nscaling},
	\begin{equation}\label{eq:proof:thm:g}
		\begin{aligned}
			g(v_1,\dots,v_n) &= \E f(U_1(v_1), \dots, U_n(v_n))\\
			&= P_{N(M(\lambda) v_1),\dots, N(M(\lambda) v_n)} (z_1, \dots, z_n)\\
			&= P_{M(\lambda) \odot_F N(v_1),\dots, M(\lambda) \odot_F N(v_n)} (z_1, \dots, z_n).
		\end{aligned}
	\end{equation}
	On the other hand,
	\begin{align*}
		&P_{M(\lambda) \odot_F N(Y(t_1)),\dots, M(\lambda) \odot_F N(Y(t_n))}  (z_1, \dots, z_n)\\
		&\hspace{4cm}= \E \left( z_1^{T_{N(Y(t_1))}(-\log M(\lambda))} \cdots z_n^{T_{N(Y(t_n))}(-\log M(\lambda)}\right)\\
		&\hspace{4cm}= \E f(U_1(Y(t_1)), \dots, U_n(Y(t_n))),
	\end{align*}
	where $f$ is as above and $U_i(u)=T_{N(u)}(-\log M(\lambda))$, $i=1,\dots,n$. Using Lemma \ref{lemma:conditioning} for $V_i=Y(t_i)$ gives
	\begin{equation*}
		P_{M(\lambda) \odot_F N(Y(t_1)),\dots, M(\lambda) \odot_F N(Y(t_n))} (z_1, \dots, z_n) = \E g(Y(t_1), \dots, Y(t_n)),
	\end{equation*}
	where $g$ is as in \eqref{eq:proof:thm:g}. We conclude that $\{ N(M(\lambda) Y(t)) \} \overset{d}{=} \{ M(\lambda) \odot_F N(Y(t)) \}$, therefore
	\begin{equation*}
		\{ X(\lambda t) \} \overset{d}{=} \{ M(\lambda) \odot_F N(Y(t)) \}.
	\end{equation*}
\end{proof}

The previous construction may seem very specific at first. However, it was shown in Ref.~\cite{vanharn1985selfsimilar} that all $\N_0$-valued self-similar processes have this form. The argument is based on the variant of Lamperti's limit theorem \cite{lamperti1962semi} which does not have its analogue for multifractal processes. See also Remark \ref{rem:directcosntr} for another construction of the process satisfying a more general definition of multifractality.

Before we give some specific examples, we establish the scaling of moments. Recall that since $Y$ is multifractal in the classical sense, the scaling property \eqref{eq:momscal} holds. The next theorem shows that the integer order moments of the process constructed in Theorem \ref{thm:constr} scale in time as a linear combination of powers of $t$.

\begin{theorem}\label{thm:mom}
	Let $X$ be the process defined in Theorem \ref{thm:constr}. For $n\in \N$ such that $\E Y(t)^n<\infty$ it holds that
	\begin{equation}\label{eq:Xconstr:mom}
		\E (X(t))^n = \sum_{k=1}^n C_n(k) t^{\tau(k)}, \quad t\in (0,T],
	\end{equation}
	where $\tau$ is the scaling function of $Y$ and $C_n(k)$ are some constants.
\end{theorem}

\begin{proof}
	For $r\in \N$, let $\E Z^{[r]}=\E Z(Z-1)\cdots (Z-r+1)$ denote the factorial moment of random variable $Z$. We start by considering factorial moments of $N(t)$. From \eqref{eq:defPP}, we have $P_{N(t)}(z)=P_{V(t)}(G(z))$. An application of Fa\`a di Bruno's formula (see e.g.~Ref.~\cite{comtet1974}) gives
	\begin{equation}
		P_{N(t)}^{(r)}(z) = \sum_{k=1}^r P_{V(t)}^{(k)} (G(z)) B_{r,k} \left(G'(z), G''(z), \dots, G^{(r-k+1)}(z)\right),
	\end{equation}
	where $B_{r,k}$ are incomplete exponential Bell polynomials (see Ref.~\cite{comtet1974}). Since $P_{V(t)}^{(k)}(z)\allowbreak = \exp \{-t (1-z)\} t^k$, we get
	\begin{equation*}
		\E (N(t))^{[r]} = P_{N(t)}^{(r)}(1) = \sum_{k=1}^r t^k B_{r,k}\left(G'(1), G''(1), \dots, G^{(n-k+1)}(1)\right).
	\end{equation*}
	From the relation between factorial and raw moments (see Ref.~\cite[Section 5.2]{daley2003}), it follows that
	\begin{equation*}
		\E (N(t))^n = \sum_{r=1}^n \left\{ \genfrac{}{}{0pt}{}{n}{r} \right\} \sum_{k=1}^r t^k B_{r,k}(G'(1), G''(1), \dots, G^{(r-k+1)}(1)),
	\end{equation*}
	where $\left\{ \genfrac{}{}{0pt}{}{n}{r} \right\}$ denotes Stirling numbers of the second kind. By denoting 
	\begin{equation*}
		K_n(k)= \sum_{r=k}^n \left\{ \genfrac{}{}{0pt}{}{n}{r} \right\} B_{r,k}(G'(1), G''(1), \dots, G^{(r-k+1)}(1)),    
	\end{equation*}
	we can write
	\begin{equation*}
		\E (N(t))^n = \sum_{k=1}^n K_n(k) t^k.
	\end{equation*}
	By Lemma \ref{lemma:conditioning} with $f(u)=u^n$, $U_1(v)=N(v)$ and $V_1=Y(t)$, it follows that
	\begin{equation*}
		\E X(t)^n = \E (N(Y(t)))^n = \sum_{k=1}^n K_n(k) \E (Y(t))^k.
	\end{equation*}
	Note that we can apply Lemma \ref{lemma:conditioning} even though $f$ is not bounded. In fact, $f$ and $N$ are monotone, and hence \eqref{eq:lemma:cond:limit} follows by the monotone convergence theorem since we have assumed $\E (Y(t))^n<\infty$. Since $Y$ is multifractal, we have $\E (Y(t))^k = c(k) t^{\tau(k)}$ and \eqref{eq:Xconstr:mom} follows by putting $C_n(k)=K_n(k)c(k)$. Symbolic computation may be used to compute the constants explicitly in specific cases (see e.g.~\cite{gao2024symbolic} for such an example).
\end{proof}

\section{Examples}

Following the construction given in the previous section, we can obtain specific examples of integer-valued multifractal processes by choosing
\begin{itemize}
	\item a nondecreasing process $Y=\{Y(t), \, t \geq0\}$ multifractal in the classical sense,
	\item a jump distribution $G$ for the unit rate compound Poisson process $N=\{N(t), \, t \geq 0\}$.
\end{itemize}

For the multifractal process $Y$ one can take a multiplicative cascade process. These can be constructed from cascade measures as follows (see Refs.~\cite{bacry2003log,muzy2002multifractal,barral2014exact} for details).

Let $\nu$ be an arbitrary infinitely divisible distribution and $\Psi$ its characteristic exponent, $\Psi(\theta) = \log \E  e^{i\theta \nu}$. Assume that $\theta_c= \sup \{ \theta \geq 0 : \E e^{\theta \nu} < \infty \} > 1$ so that the Laplace exponent $\psi (\theta) = \log \E  e^{\theta \nu}$ is finite on $[0, \theta_c)$. Furthermore, assume that $\psi(1)=0$ so that $\E e^{\nu} = 1$. Next, let $\mathcal{L}$ be an independently scattered infinitely divisible random measure on the half-plane $\mathcal{H} = \{ (u, v) : u \in \R, v \geq 0\}$ associated to $\nu$, with control measure $\mu(du,dv)=v^{-2} du dv$ (see Ref.~\cite{rajput1989spectral} for details). In particular, for every Borel set $A \subset \mathcal{H}$ such that $\mu(A) < \infty$
\begin{equation}\label{exa:cascade:lambda}
	\E \exp \left\{i \theta \mathcal{L} (A) \right\} = e^{\Psi(\theta) \mu(A)}.
\end{equation}
Fix $T > 0$ and for $t \in \R$ and $l > 0$ define sets (cones)
\begin{equation*}
	A_l(t) = \{(u,v) : v \geq l, \ -f(v)/2 < u-t \leq  f(v)/2 \},
\end{equation*}
where
\begin{equation*}
	f(v)= \begin{cases}
		v, \ v\leq T,\\
		T, \ v> T.
	\end{cases}
\end{equation*}
For $l > 0$ we can now define a random measure on $\R$ by
\begin{equation*}
	Q_l(dt) = e^{\mathcal{L} \left( A_l(t) \right)} dt.
\end{equation*}
One can show that a.s.~$Q_l$ converges weakly to a random measure $Q$, as $l\to 0$ (see Ref.~\cite{barral2014exact} for details). This limiting measure $Q$ is called the log-infinitely divisible cascade and the \textbf{cascade process} $Y=\{Y(t),\, t \geq 0\}$ is obtained by putting $Y(t) = Q([0, t])$. This process is nondecreasing and satisfies scale invariance properties \eqref{eq:mfdef} and \eqref{eq:momscal}. The scaling function $\tau$ in \eqref{eq:momscal} is determined by the Laplace exponent $\psi$ and is given by
\begin{equation*}
	\tau(q) = q - \psi(q),
\end{equation*}
while $c$ in \eqref{eq:momscal} is $c(q)=T^{-\tau(q)} \E Y(T)^q$. Different choices of the infinitely divisible measure $\nu$ lead to different cascade processes.

The jump distribution $G$ of the compound Poisson process $N$ is determined by the semigroup $F$. If we take the semigroup $F_s(z)=1-e^{-s}+e^{-s} z$, corresponding to standard multiplication operation \eqref{eq:bin-thin}, then we get from \eqref{eq:Gdef} that
\begin{equation*}
	G(z) = z.
\end{equation*}
Hence, $N$ is a standard unit rate Poisson process with unit jumps. We will fix this choice in the following examples. More examples of semigroups $F$ and corresponding $G$ can be found in Ref.~\cite{vaHarnVervaat1982}.

\medskip

If we take $Y$ to be some cascade process and $N$ a unit rate Poisson process as described, then the process $X(t)=N(Y(t))$ will be $F$-multifractal with $F_s(z)=1-e^{-s}+e^{-s} z$. Since $N$ has unit jumps, as in the proof of Theorem \ref{thm:mom} we have that $\E (N(t))^{[r]} = t^r$ and
\begin{equation*}
	\E (N(t))^n = \sum_{r=1}^n \left\{ \genfrac{}{}{0pt}{}{n}{r} \right\} t^r.
\end{equation*}
If $Y$ is a cascade process, then for $t\in (0,T]$
\begin{equation}\label{eq:mainex:momscal}
	\E (X(t))^n = \sum_{r=1}^n \left\{ \genfrac{}{}{0pt}{}{n}{r} \right\} \E (Y(t))^r = \sum_{r=1}^n \left( \left\{ \genfrac{}{}{0pt}{}{n}{r} \right\} T^{-\tau(r)} \E Y(T)^r \right) t^{\tau(r)}.
\end{equation}
However, factorial moments have a particularly simple expression. By the same argument as in the proof of Theorem \ref{thm:mom}, we get that
\begin{equation}\label{eq:mainex:facmomscal}
	\E X(t)^{[n]} = \E (N(Y(t)))^{[n]} = \E (Y(t))^n = T^{-\tau(n)} \E Y(T)^n t^{\tau(n)}.
\end{equation}

We also note that since $N$ and $Y$ both have stationary increments, so does the compound process $X$.

\begin{example}\label{ex:log-normal}
	Let $\nu$ be a normal distribution so that
	\begin{equation*}
		\psi (\theta) = \log \E  e^{\theta \nu} = m \theta + \frac{\lambda^2}{2} \theta^2, \quad \theta \in \R.
	\end{equation*}
	From the condition $\psi(1)=0$ we have that $m=-\lambda^2/2$. The resulting process $Y$ is called the log-normal multiplicative cascade process and it has the scaling function
	\begin{equation*}
		\tau(q) = q \left( 1+ \frac{\lambda^2}{2} \right) - \frac{\lambda^2}{2} q^2.
	\end{equation*}
	Taking $N$ to be a unit rate Poisson process, from Theorem \ref{thm:constr} we obtain an integer-valued multifractal process $X(t)=N(Y(t))$.
\end{example}

\begin{example}
	Let $\nu$ be a compound Poisson distribution such that
	\begin{equation*}
		\psi(\theta) = \log \E  e^{i\theta \nu} = m \theta + C \int_{\R} (e^{\theta z} - 1) \p_V(dz),
	\end{equation*}
	where $C>0$, $V$ is a random variable with distribution $\p_V$ and $m\in \R$ is such that $\psi(1)=0$. Let $W=e^V$ and assume that $\theta_c= \sup \{ \theta \geq 0 : \E W^{\theta} < \infty \} > 1$. The cascade process $Y$ obtained in this way is called the log-compound Poisson cascade process and its  scaling function is
	\begin{equation*}
		\tau(q) = q \left( 1- m \right) - C \left( \E W^q - 1\right).
	\end{equation*}
	With $N$ a unit rate Poisson process, we obtain from Theorem \ref{thm:constr} an integer-valued multifractal process.
\end{example}

\section{Simulations}
In this section, we provide empirical analysis of data simulated from time-changed models. We will restrict our attention to the process constructed from the log-normal cascade process discussed in Example \ref{ex:log-normal}.

\begin{figure}[!ht]
	\centering
	\includegraphics[width=0.7\textwidth]{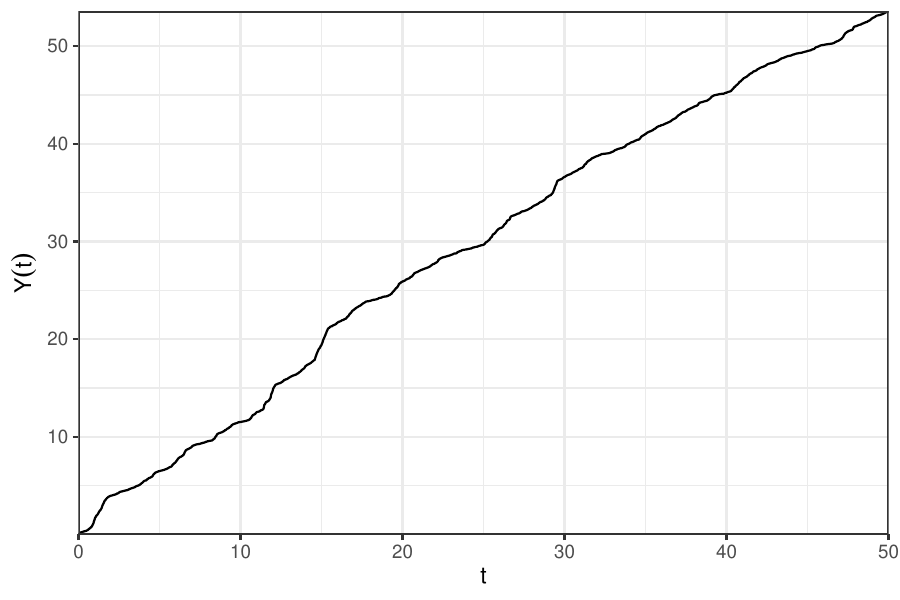}
	\caption{Simulated trajectory of the log-normal cascade process $Y$.}
	\label{fig:trajectory:cascade}
\end{figure}

\begin{figure}[!ht]
	\centering
	\includegraphics[width=0.7\textwidth]{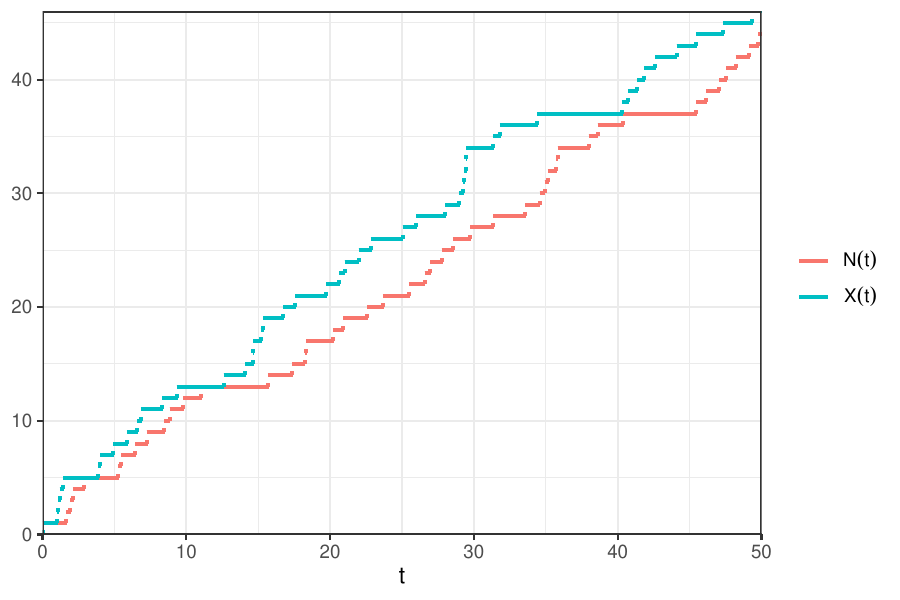}
	\caption{Simulated trajectory of the Poisson process $N$ and compound integer-valued multifractal process $X$ from Example \ref{ex:log-normal}.}
	\label{fig:trajectory:N}
\end{figure}

Figure \ref{fig:trajectory:cascade} shows the trajectory of the log-normal cascade process $Y$. For the simulation we take $\lambda^2=0.05$ and $T=50$ and use the method explained in Ref.~\cite{chainais2005} (see also Ref.~\cite{muzy2002multifractal}). Figure \ref{fig:trajectory:N} shows the trajectory of the simulated unit rate Poisson process $N$ and the process $X(t)=N(Y(t))$ obtained by changing the time of the Poisson process with the cascade process $Y$, as in Theorem \ref{thm:constr}.

We also investigate the behavior of moments by estimating the scaling function using empirical moments. Let $U=\{U(t),\, t \in [0,T]\}$ denote a process with stationary increments. By dividing the interval $[0,T]$ into $\lfloor T / t \rfloor$ blocks of length $t$, the moment $\E |U(t)|^q$ can be estimated with
\begin{equation*}
	S_q(T,t) = \frac{1}{\lfloor T /  t \rfloor} \sum_{i=1}^{\lfloor T /  t \rfloor} \left| U ( i  t) -  U ( (i-1)  t) \right|^q.
\end{equation*}
The function $S_q$ is usually called the partition function or the empirical structure function. If the moments of the process $U$ scale as in \eqref{eq:momscal}, then we expect that, for fixed $q$ and some $0\leq {t}_1 < \cdots < {t}_N \leq T$, the points $(\log t_i , \log S_q(T, t_i))$, $i=1,\dots,N$, exhibit a linear pattern. In this case, we can estimate $\tau(q)$ as the slope of the simple linear regression of $\ln S_q(T, t)$ on $\ln {t}$ using the well-known formula for the slope of the linear regression line:
\begin{equation}\label{eq:tauhat}
	\hat{\tau}_{N,T}(q) = \frac{\sum_{i=1}^{N}  \ln {t_i}  \ln S_q(n,t_i) - \frac{1}{N} \sum_{i=1}^{N} \ln {t_i} \sum_{j=1}^{N} \ln S_q(n,t_i) }{ \sum_{i=1}^{N} \left(\ln {t_i}\right)^2 - \frac{1}{N} \left( \sum_{i=1}^{N} \ln {t_i} \right)^2 },
\end{equation}
where $N$ is the number of time points chosen in the regression. This can be repeated for a range of values of $q$ to obtain the plot of the estimated scaling function. For more details on this methodology see Refs.~\cite{FCM1997multifractalityDEM,grahovac2014detecting} and the references therein. Alternatively, wavelet methods could be used to analyze scale invariance (see Refs.~\cite{abry2015irregularities,jaffard2001wavelets} for the survey and \cite{grahovac2018bounds} for other methods).

This approach works well with multifractal processes satisfying \eqref{eq:momscal}. However, the moments of the integer-valued process $X$ scale as a linear combination of powers of $t$. The coefficients in \eqref{eq:mainex:momscal} will determine which of these powers will dominate and these depend on how large $t$ and $T$ are.

Since the factorial moments of $X$ exhibit exact power law scaling \eqref{eq:mainex:facmomscal}, we can try to estimate the scaling function from estimated factorial moments. To this end, we define for $n\in \N$ the factorial moments partition function
\begin{equation}\label{eq:partitionfun:fac}
	S_{[n]}(T,t) = \frac{1}{\lfloor T /  t \rfloor} \sum_{i=1}^{\lfloor T /  t \rfloor} m^{[n]} \left( X( i  t) - X((i-1)t)\right),
\end{equation}
where
\begin{equation*}
	m^{[n]} (x) = x(x-1)\cdots (x-n+1).
\end{equation*}
Similarly as for the sample absolute moments, one can estimate the scaling function as in \eqref{eq:tauhat} from the points $(\log t_i , \allowbreak \log S_{[n]}(T, t_i))$, $i=1,\dots,N$.

\begin{figure}[!ht]
	\centering
	\includegraphics[width=0.7\textwidth]{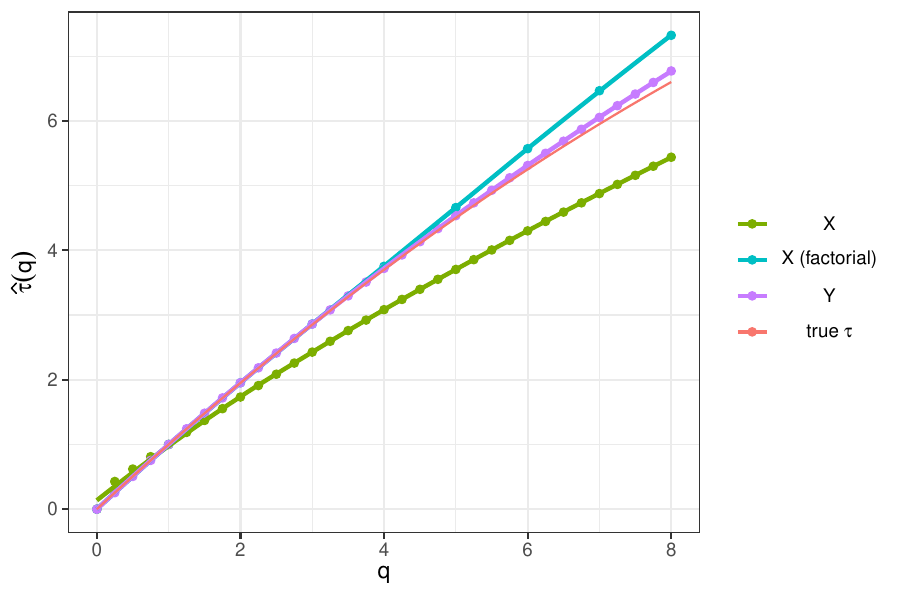}
	\caption{True and estimated scaling functions ($T=10000$).}
	\label{fig:sf:longer}
\end{figure}

\begin{figure}[!ht]
	\centering
	\includegraphics[width=0.7\textwidth]{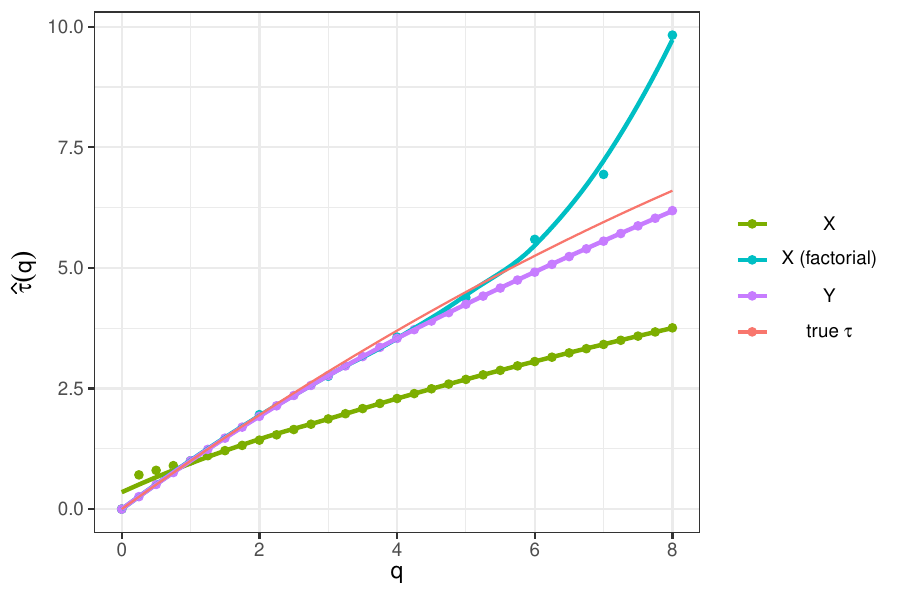}
	\caption{True and estimated scaling functions $T=100$.}
	\label{fig:sf:shorter}
\end{figure}

Figures \ref{fig:sf:longer} and \ref{fig:sf:shorter} compare the true scaling function with estimates for the cascade process $Y$ and the compound process $X$, using both sample moments and sample factorial moments. The cascade process is log-normal with parameter $\lambda^2 = 0.05$, simulated with time step $0.01$. In Figure \ref{fig:sf:longer} the path length is $T = 10000$, while in Figure \ref{fig:sf:shorter} it is $T = 100$.

For $Y$, the estimated scaling function closely matches the true curve in both cases. In contrast, applying the same moment-based estimation to $X$ shows noticeable discrepancies. This behavior is expected: the moments of $X$ do not follow an exact power-law, and the deviations depend on the coefficients in \eqref{eq:mainex:momscal}, which themselves vary with $T$.

When factorial moments are used instead, the estimated scaling function aligns much better with the true form, in line with \eqref{eq:mainex:facmomscal}. However, higher-order factorial moments require evaluating large factorials, which can introduce numerical instabilities. Moreover, \eqref{eq:partitionfun:fac} may become zero for large $n$, hence only a fraction of points is used to estimate the scaling function for these $n$. These instabilities are particularly visible in Figure \ref{fig:sf:shorter}. In practical applications, scaling functions for integer-valued data can be estimated using either moments or factorial moments. The parameters of the proposed model may then be obtained by fitting the estimated scaling function to its theoretical form (see Ref.~\cite{grahovac2014detecting} for details on this methodology and illustrative examples).

\section{Conclusion}

In this paper we introduced a multifractal framework for integer-valued processes, extending multifractal models beyond the classical real-valued setting. By employing the thinning operation as a discrete analogue of scalar multiplication, we formalized a definition of multifractality for count data. The construction via multifractal time-changes of compound Poisson processes provides a tractable and versatile class of examples, for which moment scaling laws were obtained and numerical illustrations carried out.  

Although count data arise in many contexts, their multifractal characteristics have not been systematically explored. Moreover, in applications to count data, multifractal methods have typically relied on real-valued models, which serve only as approximations. The results presented here provide a foundation for studying count data that display variability across multiple scales, opening the way to applications in a variety of fields.

The proposed framework leaves several questions open. In particular, the problem of estimating the scaling function from count data is highly relevant, and a more detailed analysis of asymptotic properties is needed. It is also natural to investigate extensions of these ideas to multivariate and spatial processes, as well as to random measures. Finally, applications to empirical data sets from epidemiology, finance, and physics could provide valuable validation of the models and suggest possible refinements.

\bigskip
\bigskip

%
%
%

\bibliographystyle{mcaomplain} 
\bibliography{References}

\providecommand{\bysame}{\leavevmode\hbox to3em{\hrulefill}\thinspace}
\providecommand{\noopsort}[1]{}
\providecommand{\arxiv}[1]{\href{http://www.arxiv.org/abs/#1}{arXiv~#1}}
\providecommand{\doi}[1]{\url{https://doi.org/#1}}
\providecommand{\href}[2]{#2}
\bibliofont{\begin{thebibliography}{10}

\bibitem{abry2015irregularities}
\bgroup\scshape{}P.~Abry\egroup{}, \bgroup\scshape{}S.~Jaffard\egroup{}, and
  \bgroup\scshape{}H.~Wendt\egroup{},
  \href{https://doi.org/10.1142/9789814366076\_0003}{Irregularities and scaling
  in signal and image processing: multifractal analysis},  in \emph{Benoit
  Mandelbrot: a life in many dimensions}, World Scientific, 2015, pp.~31--116.

\bibitem{anh2009multifractal}
\bgroup\scshape{}V.~V. Anh\egroup{}, \bgroup\scshape{}N.~N. Leonenko\egroup{},
  and \bgroup\scshape{}N.-R. Shieh\egroup{},
  \href{https://doi.org/10.1080/07362990802679091}{Multifractal products of
  stationary diffusion processes},  \emph{Stoch. Anal. Appl.} \textbf{27}
  (2009), no.~3, 475--499.

\bibitem{apolinario2022dynamical}
\bgroup\scshape{}G.~B. Apolin{\'a}rio\egroup{},
  \bgroup\scshape{}L.~Chevillard\egroup{}, and \bgroup\scshape{}J.-C.
  Mourrat\egroup{}, \href{https://doi.org/10.1007/s10955-021-02867-2}{Dynamical
  fractional and multifractal fields},  \emph{Journal of Statistical Physics}
  \textbf{186} (2022), no.~1, 15.

\bibitem{asmussenhering1983}
\bgroup\scshape{}S.~Asmussen\egroup{} and \bgroup\scshape{}H.~Hering\egroup{},
  \emph{Branching processes}, \emph{Progress in Probability and Statistics}
  \textbf{3}, Birkh\"{a}user Boston, Inc., Boston, MA, 1983.

\bibitem{AthreyaNey1972}
\bgroup\scshape{}K.~B. Athreya\egroup{} and \bgroup\scshape{}P.~E.
  Ney\egroup{}, \emph{Branching processes}, \emph{Die Grundlehren der
  mathematischen Wissenschaften} \textbf{196}, Springer-Verlag, New
  York-Heidelberg, 1972.

\bibitem{bacry2008continuous}
\bgroup\scshape{}E.~Bacry\egroup{}, \bgroup\scshape{}A.~Kozhemyak\egroup{}, and
  \bgroup\scshape{}J.-F. Muzy\egroup{},
  \href{https://doi.org/10.1016/j.jedc.2007.01.024}{Continuous cascade models
  for asset returns},  \emph{Journal of Economic Dynamics and Control}
  \textbf{32} (2008), no.~1, 156--199.

\bibitem{bacry2003log}
\bgroup\scshape{}E.~Bacry\egroup{} and \bgroup\scshape{}J.~F. Muzy\egroup{},
  \href{https://doi.org/10.1007/s00220-003-0827-3}{Log-infinitely divisible
  multifractal processes},  \emph{Communications in Mathematical Physics}
  \textbf{236} (2003), no.~3, 449--475.

\bibitem{BARANOWSKI2019318}
\bgroup\scshape{}P.~Baranowski\egroup{}, \bgroup\scshape{}M.~Gos\egroup{},
  \bgroup\scshape{}J.~Krzyszczak\egroup{}, \bgroup\scshape{}K.~Siwek\egroup{},
  \bgroup\scshape{}A.~Kieliszek\egroup{}, and
  \bgroup\scshape{}P.~Tkaczyk\egroup{},
  \href{https://doi.org/10.1016/j.chaos.2019.07.008}{Multifractality of
  meteorological time series for {P}oland on the base of {MERRA}-2 data},
  \emph{Chaos, Solitons \& Fractals} \textbf{127} (2019), 318--333.

\bibitem{barral2014exact}
\bgroup\scshape{}J.~Barral\egroup{} and \bgroup\scshape{}X.~Jin\egroup{},
  \href{https://doi.org/10.1007/s00440-013-0534-8}{On exact scaling
  log-infinitely divisible cascades},  \emph{Probability Theory and Related
  Fields} \textbf{160} (2014), no.~3-4, 521--565.

\bibitem{barral2002multifractal}
\bgroup\scshape{}J.~Barral\egroup{} and \bgroup\scshape{}B.~B.
  Mandelbrot\egroup{},
  \href{https://doi.org/10.1007/s004400200220}{Multifractal products of
  cylindrical pulses},  \emph{Probability Theory and Related Fields}
  \textbf{124} (2002), no.~3, 409--430.

\bibitem{baile2021}
\bgroup\scshape{}R.~Baïle\egroup{}, \bgroup\scshape{}J.-F. Muzy\egroup{}, and
  \bgroup\scshape{}X.~Silvani\egroup{},
  \href{https://doi.org/https://doi.org/10.1016/j.physa.2020.125697}{Multifractal
  point processes and the spatial distribution of wildfires in {F}rench
  {M}editerranean regions},  \emph{Physica A: Statistical Mechanics and its
  Applications} \textbf{568} (2021), 125697.

\bibitem{boucher2013time}
\bgroup\scshape{}J.-P. Boucher\egroup{} and
  \bgroup\scshape{}D.~Hainaut\egroup{}, Time series of count data using
  multifractal process,  \emph{Preprint. (UQAM)} (2013). Available at
  \url{http://archipel.uqam.ca/id/eprint/6967}.

\bibitem{boucher2015time}
\bgroup\scshape{}J.-P. Boucher\egroup{} and
  \bgroup\scshape{}D.~Hainaut\egroup{}, Time series of correlated count data
  using multifractal process,  \emph{Preprint. (UQAM)} (2015). Available at
  \url{http://archipel.uqam.ca/id/eprint/6968}.

\bibitem{cameron_trivedi_1998}
\bgroup\scshape{}A.~C. Cameron\egroup{} and \bgroup\scshape{}P.~K.
  Trivedi\egroup{}, \emph{Regression analysis of count data}, second ed.,
  \emph{Econometric Society Monographs} \textbf{53}, Cambridge University
  Press, Cambridge, 2013.

\bibitem{chainais2005}
\bgroup\scshape{}P.~Chainais\egroup{}, \bgroup\scshape{}R.~Riedi\egroup{}, and
  \bgroup\scshape{}P.~Abry\egroup{},
  \href{https://doi.org/10.1109/TIT.2004.842570}{On non-scale-invariant
  infinitely divisible cascades},  \emph{IEEE Trans. Inform. Theory}
  \textbf{51} (2005), no.~3, 1063--1083.

\bibitem{comtet1974}
\bgroup\scshape{}L.~Comtet\egroup{}, \emph{Advanced combinatorics}, D. Reidel
  Publishing Co., Dordrecht, 1974.

\bibitem{daley2003}
\bgroup\scshape{}D.~J. Daley\egroup{} and
  \bgroup\scshape{}D.~Vere-Jones\egroup{}, \emph{An introduction to the theory
  of point processes. {V}ol. {I}}, second ed., \emph{Probability and its
  Applications}, Springer-Verlag, New York, 2003.

\bibitem{davis1998}
\bgroup\scshape{}R.~A. Davis\egroup{}, \bgroup\scshape{}Y.~Wang\egroup{}, and
  \bgroup\scshape{}W.~T.~M. Dunsmuir\egroup{}, Modeling time series of count
  data,  in \emph{Asymptotics, nonparametrics, and time series}, \emph{Statist.
  Textbooks Monogr.} \textbf{158}, Dekker, New York, 1999, pp.~63--113.

\bibitem{duchon2012forecasting}
\bgroup\scshape{}J.~Duchon\egroup{}, \bgroup\scshape{}R.~Robert\egroup{}, and
  \bgroup\scshape{}V.~Vargas\egroup{},
  \href{https://doi.org/10.1111/j.1467-9965.2010.00458.x}{Forecasting
  volatility with the multifractal random walk model},  \emph{Math. Finance}
  \textbf{22} (2012), no.~1, 83--108.

\bibitem{ebrahimkhanlou2016multifractal}
\bgroup\scshape{}A.~Ebrahimkhanlou\egroup{},
  \bgroup\scshape{}A.~Farhidzadeh\egroup{}, and
  \bgroup\scshape{}S.~Salamone\egroup{},
  \href{https://doi.org/10.1177/1475921715624502}{Multifractal analysis of
  crack patterns in reinforced concrete shear walls},  \emph{Structural Health
  Monitoring} \textbf{15} (2016), no.~1, 81--92.

\bibitem{feng2025bilinear}
\bgroup\scshape{}C.-H. Feng\egroup{}, \bgroup\scshape{}B.~Tian\egroup{}, and
  \bgroup\scshape{}X.-T. Gao\egroup{},
  \href{https://doi.org/10.1007/s12346-025-01241-x}{Bilinear {B}{\"a}cklund
  transformations, as well as {$N$}-soliton, breather, fission/fusion and
  hybrid solutions for a $(3+ 1)$-dimensional integrable wave equation in a
  fluid},  \emph{Qualitative Theory of Dynamical Systems} \textbf{24} (2025),
  100.

\bibitem{FCM1997multifractalityDEM}
\bgroup\scshape{}A.~Fisher\egroup{}, \bgroup\scshape{}L.~Calvet\egroup{}, and
  \bgroup\scshape{}B.~B. Mandelbrot\egroup{}, Multifractality of
  {D}eutschemark/{US} {D}ollar exchange rates,  \emph{Cowles Foundation
  discussion paper} \textbf{1166} (1997). Available at
  \url{https://ssrn.com/abstract=78628}.

\bibitem{frisch1985fully}
\bgroup\scshape{}U.~Frisch\egroup{} and \bgroup\scshape{}G.~Parisi\egroup{}, On
  the singularity structure of fully developed turbulence,  in \emph{Turbulence
  and predictability in geophysical fluid dynamics and climate dynamics,
  Proceed. Intern. School of Physics, Varenna, Italy}, 1985, pp.~84--87.

\bibitem{gao2024symbolic}
\bgroup\scshape{}X.-Y. Gao\egroup{}, Open-ocean shallow-water dynamics via a
  {$(2+1)$}-dimensional generalized variable-coefficient {Hirota-Satsuma-Ito}
  system: Oceanic auto-{Bä}cklund transformation and oceanic solitons,
  \emph{China Ocean Engineering} \textbf{39} (2025), no.~3, 541–547.

\bibitem{gao2025inhomogeneity}
\bgroup\scshape{}X.-Y. Gao\egroup{}, \bgroup\scshape{}J.-G. Liu\egroup{}, and
  \bgroup\scshape{}G.-W. Wang\egroup{},
  \href{https://doi.org/10.1016/j.aml.2025.109615}{Inhomogeneity, magnetic
  auto-{B}{\"a}cklund transformations and magnetic solitons for a generalized
  variable-coefficient {Kraenkel-Manna-Merle} system in a deformed ferrite},
  \emph{Applied Mathematics Letters} \textbf{171} (2025), 109615.

\bibitem{grahovac2020}
\bgroup\scshape{}D.~Grahovac\egroup{},
  \href{https://doi.org/10.1016/j.chaos.2020.109735}{Multifractal processes:
  definition, properties and new examples},  \emph{Chaos, Solitons \& Fractals}
  \textbf{134} (2020), 109735, 11.

\bibitem{grahovac2014detecting}
\bgroup\scshape{}D.~Grahovac\egroup{} and \bgroup\scshape{}N.~N.
  Leonenko\egroup{},
  \href{https://doi.org/10.1016/j.chaos.2014.04.016}{Detecting multifractal
  stochastic processes under heavy-tailed effects},  \emph{Chaos, Solitons \&
  Fractals} \textbf{65} (2014), 78--89.

\bibitem{grahovac2018bounds}
\bgroup\scshape{}D.~Grahovac\egroup{} and \bgroup\scshape{}N.~N.
  Leonenko\egroup{}, \href{https://doi.org/10.1142/S0218348X1850055X}{Bounds on
  the support of the multifractal spectrum of stochastic processes},
  \emph{Fractals} \textbf{26} (2018), no.~04, 1850055.

\bibitem{vanharn1985selfsimilar}
\bgroup\scshape{}K.~van Harn\egroup{} and \bgroup\scshape{}F.~W.
  Steutel\egroup{},
  \href{https://doi.org/10.1080/15326348508807010}{Integer-valued self-similar
  processes},  \emph{Comm. Statist. Stochastic Models} \textbf{1} (1985),
  no.~2, 191--208.

\bibitem{vaHarnVervaat1982}
\bgroup\scshape{}K.~van Harn\egroup{}, \bgroup\scshape{}F.~W. Steutel\egroup{},
  and \bgroup\scshape{}W.~Vervaat\egroup{},
  \href{https://doi.org/10.1007/BF00537228}{Self-decomposable discrete
  distributions and branching processes},  \emph{Z. Wahrsch. Verw. Gebiete}
  \textbf{61} (1982), no.~1, 97--118.

\bibitem{harris1963}
\bgroup\scshape{}T.~E. Harris\egroup{}, \emph{The theory of branching
  processes}, \emph{Die Grundlehren der Mathematischen Wissenschaften}
  \textbf{119}, Springer-Verlag, Berlin, 1963.

\bibitem{Hilbe_2014}
\bgroup\scshape{}J.~M. Hilbe\egroup{}, \emph{Modeling count data}, Cambridge
  University Press, Cambridge, 2014.

\bibitem{jaffard2001wavelets}
\bgroup\scshape{}S.~Jaffard\egroup{}, \bgroup\scshape{}Y.~Meyer\egroup{}, and
  \bgroup\scshape{}R.~D. Ryan\egroup{}, \emph{Wavelets: tools for science and
  technology}, SIAM, Philadelphia, 2001.

\bibitem{jiang2019multifractal}
\bgroup\scshape{}Z.-Q. Jiang\egroup{}, \bgroup\scshape{}W.-J. Xie\egroup{},
  \bgroup\scshape{}W.-X. Zhou\egroup{}, and
  \bgroup\scshape{}D.~Sornette\egroup{},
  \href{https://doi.org/10.1088/1361-6633/ab42fb}{Multifractal analysis of
  financial markets: a review},  \emph{Rep. Progr. Phys.} \textbf{82} (2019),
  no.~12, 125901, 105.

\bibitem{joseph2021multifractal}
\bgroup\scshape{}A.~J. Joseph\egroup{} and \bgroup\scshape{}P.~N.
  Pournami\egroup{},
  \href{https://doi.org/10.1016/j.chaos.2021.111301}{Multifractal theory based
  breast tissue characterization for early detection of breast cancer},
  \emph{Chaos, Solitons \& Fractals} \textbf{152} (2021), 111301.

\bibitem{kalamaras2017multifractal}
\bgroup\scshape{}N.~Kalamaras\egroup{},
  \bgroup\scshape{}K.~Philippopoulos\egroup{},
  \bgroup\scshape{}D.~Deligiorgi\egroup{}, \bgroup\scshape{}C.~Tzanis\egroup{},
  and \bgroup\scshape{}G.~Karvounis\egroup{},
  \href{https://doi.org/10.1016/j.chaos.2017.03.003}{Multifractal scaling
  properties of daily air temperature time series},  \emph{Chaos, Solitons \&
  Fractals} \textbf{98} (2017), 38--43.

\bibitem{laib2018multifractal}
\bgroup\scshape{}M.~Laib\egroup{}, \bgroup\scshape{}J.~Golay\egroup{},
  \bgroup\scshape{}L.~Telesca\egroup{}, and
  \bgroup\scshape{}M.~Kanevski\egroup{},
  \href{https://doi.org/10.1016/j.chaos.2018.02.024}{Multifractal analysis of
  the time series of daily means of wind speed in complex regions},
  \emph{Chaos, Solitons \& Fractals} \textbf{109} (2018), 118--127.

\bibitem{lamperti1962semi}
\bgroup\scshape{}J.~Lamperti\egroup{},
  \href{https://doi.org/10.2307/1993933}{Semi-stable stochastic processes},
  \emph{Transactions of the American Mathematical Society} \textbf{104} (1962),
  62--78.

\bibitem{leonenko2007}
\bgroup\scshape{}N.~N. Leonenko\egroup{}, \bgroup\scshape{}V.~Savani\egroup{},
  and \bgroup\scshape{}A.~A. Zhigljavsky\egroup{}, Autoregressive negative
  binomial processes,  \emph{Ann. I.S.U.P.} \textbf{51} (2007), no.~1-2,
  25--47.

\bibitem{liu2025n}
\bgroup\scshape{}H.-D. Liu\egroup{}, \bgroup\scshape{}B.~Tian\egroup{},
  \bgroup\scshape{}Y.-Q. Chen\egroup{}, \bgroup\scshape{}C.-D. Cheng\egroup{},
  and \bgroup\scshape{}X.-T. Gao\egroup{},
  \href{https://doi.org/10.1007/s11071-024-10397-1}{$n$-soliton, {$H$}th-order
  breather, hybrid and multi-pole solutions for a generalized
  variable-coefficient {G}ardner equation with an external force in a plasma or
  fluid},  \emph{Nonlinear Dynamics} \textbf{113} (2025), no.~4, 3655--3672.

\bibitem{lopes2009fractal}
\bgroup\scshape{}R.~Lopes\egroup{} and \bgroup\scshape{}N.~Betrouni\egroup{},
  \href{https://doi.org/10.1016/j.media.2009.05.003}{Fractal and multifractal
  analysis: a review},  \emph{Medical Image Analysis} \textbf{13} (2009),
  no.~4, 634--649.

\bibitem{lovejoy2013weather}
\bgroup\scshape{}S.~Lovejoy\egroup{} and
  \bgroup\scshape{}D.~Schertzer\egroup{}, \emph{The weather and climate:
  emergent laws and multifractal cascades}, Cambridge University Press,
  Cambridge, 2013.

\bibitem{mandelbrot1972}
\bgroup\scshape{}B.~B. Mandelbrot\egroup{}, Possible refinement of the
  lognormal hypothesis concerning the distribution of energy dissipation in
  intermittent turbulence,  in \emph{Statistical models and turbulence}
  (\bgroup\scshape{}M.~Rosenblatt\egroup{} and
  \bgroup\scshape{}C.~Van~Atta\egroup{}, eds.), \emph{Lecture Notes in Physics}
  \textbf{12}, Springer, Berlin, 1972, pp.~333--351.

\bibitem{mandelbrot1997fractal}
\bgroup\scshape{}B.~B. Mandelbrot\egroup{}, \emph{The fractal geometry of
  nature}, W. H. Freeman and Co., San Francisco, CA, 1982.

\bibitem{mandelbrot1997mmar}
\bgroup\scshape{}B.~B. Mandelbrot\egroup{},
  \bgroup\scshape{}A.~Fisher\egroup{}, and \bgroup\scshape{}L.~Calvet\egroup{},
  A multifractal model of asset returns,  \emph{Cowles Foundation discussion
  paper} \textbf{1164} (1997). Available at
  \url{https://ssrn.com/abstract=78588}.

\bibitem{martsepp2022dependence}
\bgroup\scshape{}M.~Martsepp\egroup{}, \bgroup\scshape{}T.~Laas\egroup{},
  \bgroup\scshape{}K.~Laas\egroup{}, \bgroup\scshape{}J.~Priimets\egroup{},
  \bgroup\scshape{}S.~Tõkke\egroup{}, and \bgroup\scshape{}V.~Mikli\egroup{},
  \href{https://doi.org/10.1016/j.chaos.2022.111811}{Dependence of multifractal
  analysis parameters on the darkness of a processed image},  \emph{Chaos,
  Solitons \& Fractals} \textbf{156} (2022), 111811.

\bibitem{mckenzie2003}
\bgroup\scshape{}E.~McKenzie\egroup{},
  \href{https://doi.org/10.1016/S0169-7161(03)21018-X}{Discrete variate time
  series},  in \emph{Stochastic processes: modelling and simulation},
  \emph{Handbook of Statistics} \textbf{21}, North-Holland, Amsterdam, 2003,
  pp.~573--606.

\bibitem{gadre2003multifractal}
\bgroup\scshape{}K.~P. Murali\egroup{}, \bgroup\scshape{}V.~M. Gadre\egroup{},
  and \bgroup\scshape{}U.~B. Desai\egroup{}, \emph{Multifractal based network
  traffic modeling}, Springer New York, NY, 2003.

\bibitem{muzy2002multifractal}
\bgroup\scshape{}J.-F. Muzy\egroup{} and \bgroup\scshape{}E.~Bacry\egroup{},
  \href{https://doi.org/10.1103/PhysRevE.66.056121}{Multifractal stationary
  random measures and multifractal random walks with log-infinitely divisible
  scaling laws},  \emph{Physical Review E} \textbf{66} (2002), no.~5, 056121.

\bibitem{pavlov2016multifractality}
\bgroup\scshape{}A.~Pavlov\egroup{},
  \bgroup\scshape{}O.~Semyachkina-Glushkovskaya\egroup{},
  \bgroup\scshape{}O.~Pavlova\egroup{},
  \bgroup\scshape{}A.~Abdurashitov\egroup{},
  \bgroup\scshape{}G.~Shihalov\egroup{}, \bgroup\scshape{}E.~Rybalova\egroup{},
  and \bgroup\scshape{}S.~Sindeev\egroup{},
  \href{https://doi.org/10.1016/j.chaos.2016.06.002}{Multifractality in
  cerebrovascular dynamics: an approach for mechanisms-related analysis},
  \emph{Chaos, Solitons \& Fractals} \textbf{91} (2016), 210--213.

\bibitem{rajput1989spectral}
\bgroup\scshape{}B.~S. Rajput\egroup{} and
  \bgroup\scshape{}J.~Rosinski\egroup{},
  \href{https://doi.org/10.1007/BF00339998}{Spectral representations of
  infinitely divisible processes},  \emph{Probability Theory and Related
  Fields} \textbf{82} (1989), no.~3, 451--487.

\bibitem{riedi2003multifractal}
\bgroup\scshape{}R.~H. Riedi\egroup{}, Multifractal processes,  in \emph{Theory
  and applications of long-range dependence}
  (\bgroup\scshape{}P.~Doukhan\egroup{},
  \bgroup\scshape{}G.~Oppenheim\egroup{}, and \bgroup\scshape{}M.~S.
  Taqqu\egroup{}, eds.), Birkh\"auser Basel, 2003, pp.~625--716.

\bibitem{robert2008hydrodynamic}
\bgroup\scshape{}R.~Robert\egroup{} and \bgroup\scshape{}V.~Vargas\egroup{},
  \href{https://doi.org/10.1007/s00220-008-0642-y}{Hydrodynamic turbulence and
  intermittent random fields},  \emph{Comm. Math. Phys.} \textbf{284} (2008),
  no.~3, 649--673.

\bibitem{sato1999levy}
\bgroup\scshape{}K.~Sato\egroup{}, \emph{L\'{e}vy processes and infinitely
  divisible distributions}, \emph{Cambridge Studies in Advanced Mathematics}
  \textbf{68}, Cambridge University Press, Cambridge, 1999.

\bibitem{steutel1979}
\bgroup\scshape{}F.~W. Steutel\egroup{} and \bgroup\scshape{}K.~van
  Harn\egroup{}, Discrete analogues of self-decomposability and stability,
  \emph{Ann. Probab.} \textbf{7} (1979), no.~5, 893--899.

\bibitem{steutelvanharn2004book}
\bgroup\scshape{}F.~W. Steutel\egroup{} and \bgroup\scshape{}K.~van
  Harn\egroup{}, \emph{Infinite divisibility of probability distributions on
  the real line}, \emph{Monographs and Textbooks in Pure and Applied
  Mathematics} \textbf{259}, Marcel Dekker, Inc., New York, 2004.

\bibitem{strogatz2024nonlinear}
\bgroup\scshape{}S.~H. Strogatz\egroup{}, \emph{Nonlinear dynamics and chaos:
  with applications to physics, biology, chemistry, and engineering}, 2 ed.,
  CRC Press, Boca Raton, 2015.

\bibitem{vehel2012fractals}
\bgroup\scshape{}J.~L. V{\'e}hel\egroup{}, \bgroup\scshape{}E.~Lutton\egroup{},
  and \bgroup\scshape{}C.~Tricot\egroup{}, \emph{Fractals in engineering: from
  theory to industrial applications}, Springer Science \& Business Media,
  London, 2012.

\bibitem{WANG2026109720}
\bgroup\scshape{}G.~Wang\egroup{}, \bgroup\scshape{}Z.~Tan\egroup{},
  \bgroup\scshape{}X.-Y. Gao\egroup{}, and \bgroup\scshape{}J.-G. Liu\egroup{},
  \href{https://doi.org/10.1016/j.aml.2025.109720}{A new $(2+1)$-dimensional
  {like-Harry-Dym} equation with derivation and soliton solutions},
  \emph{Applied Mathematics Letters} \textbf{172} (2026), 109720.

\bibitem{weiss2008}
\bgroup\scshape{}C.~H. Wei\ss\egroup{},
  \href{https://doi.org/10.1007/s10182-008-0072-3}{Thinning operations for
  modeling time series of counts---a survey},  \emph{AStA Adv. Stat. Anal.}
  \textbf{92} (2008), no.~3, 319--341.

\bibitem{wu2015low}
\bgroup\scshape{}Z.~Wu\egroup{}, \bgroup\scshape{}L.~Zhang\egroup{}, and
  \bgroup\scshape{}M.~Yue\egroup{},
  \href{https://doi.org/10.1109/TDSC.2015.2443807}{Low-rate {DoS} attacks
  detection based on network multifractal},  \emph{IEEE Transactions on
  Dependable and Secure Computing} \textbf{13} (2016), no.~5, 559--567.

\bibitem{ZHOU2012147}
\bgroup\scshape{}W.-X. Zhou\egroup{},
  \href{https://doi.org/10.1016/j.chaos.2011.11.004}{Finite-size effect and the
  components of multifractality in financial volatility},  \emph{Chaos,
  Solitons \& Fractals} \textbf{45} (2012), no.~2, 147--155.

\end{thebibliography}



%
%
%
\end{document}